\documentclass[10pt]{amsart}

\usepackage{tikz}
\usepackage{amssymb,amsmath,latexsym,graphicx}
\usepackage{charter,eucal}
\usepackage{amssymb}
\usepackage[all,cmtip]{xy}
\usepackage{rotating,multicol}
\usetikzlibrary{decorations.markings}

\newtheorem{theorem}{Theorem }[section]

\newtheorem{lemma}[theorem]{Lemma}
\newtheorem{observation}[theorem]{Observation}

\newtheorem{remark}[theorem]{Remark}
\newtheorem{corollary}[theorem]{Corollary}
\newtheorem{proposition}[theorem]{Proposition}
\newtheorem{principle}[theorem]{\textsc{Principle}}

\newcommand{\bt}{\begin{theorem}}
\newcommand{\et}{\end{theorem}}
\newcommand{\bmt}{\begin{maintheorem}}
\newcommand{\emt}{\end{maintheorem}}
\newcommand{\bc}{\begin{corollary}}
\newcommand{\bl}{\begin{lemma}}
\newcommand{\ec}{\end{corollary}}
\newcommand{\el}{\end{lemma}}
\newcommand{\bo}{\begin{observation}}
\newcommand{\eo}{\end{observation}}
\newcommand{\bp}{\begin{proposition}}
\newcommand{\ep}{\end{proposition}}
\newcommand{\br}{\begin{remark}}
\newcommand{\er}{\end{remark}}
\newcommand{\bpr}{\begin{principle}}
\newcommand{\epr}{\end{principle}}

\def\PG{\mathbf{PG}}

\def\AG{\mathbf{AG}}

\def\eop{\hspace*{\fill}$\blacksquare$}

\def\O{\mathbf{O}}

\usepackage{amssymb} 
\usepackage{amsmath} 
\usepackage[utf8]{inputenc} 
\usepackage{latexsym}

\title {\bf Generalized quadrangles of order $(s, s^2)$. IV \\ Translations, Moufang and Fong-Seitz}
\author {J. A. THAS \\ Ghent University}
\address{Ghent University, Department of Mathematics, Krijgslaan 281, S25, B-9000 Ghent, Belgium}
\email{thas.joseph@gmail.com}
\begin {document}
\maketitle
\begin{abstract}
In the period 1994-1999 Thas wrote a series of three papers on generalized quadrangles of order $(s, s^2)$. In this Part IV we classify all finite translation generalized quadrangles of order $(s, s^2)$ having a kernel of size at least 3, containing a regular line not incident with the translation point. There are several applications on generalized quadrangles of order $(s, s^2)$ having at least 2 translation points, on Moufang quadrangles, and concerning the theorem of Fong and Seitz classifying all groups with a $BN$-pair of rank 2.

\end{abstract}
\section{Introduction}
In 1976 Tits \cite {T:76} shows that each finite Moufang generalized quadrangle $\mathcal {S}$ is classical or dual classical. This follows from the classification of all groups with a $BN$-pair of rank 2 by Fong and Seitz \cite {F&S:73,F&S:74}. In their monograph "Finite Generalized Quadrangles", the first edition of which appeared in 1984, Payne and Thas \cite {P&JAT:84,P&JAT:09} did an almost successful attempt to prove the Moufang theorem in a pure geometric way, using just some elementary group theory. As a corollary there would follow a pure geometric proof of the theorem of Fong and Seitz, namely the case where the Weyl group is $D_8$.
\par The case that was missing in the work of Payne and Thas is when the generalized quadrangle has order $(s, s^2), s \not = 1$, and each point is a translation point. In 1991 Kantor \cite {K:91} succeeded in handling this and compared to Fong and Seitz the amount of group theory employed was very small: information was required concerning a relatively restricted type of 2-transitive permutation group, together with a fact about very small degree representations of such a group. In \cite {KT:01} K. Thas gave an alternative proof relying on the classification of finite groups with a split $BN$-pair of rank 1, a classification independently obtained in 1972 by Hering, Kantor and Seitz \cite {HKS:72} and Shult \cite {S:72}.
\par In this paper we will give a pure geometric proof, without any group theory. In the even case we exclude the case where the "kernel" has just 2 elements; we hope to be able to solve also that case. In fact we obtain much more general results, characterizing classes of generalized quadrangles which are not necessarily classical nor dual classical, and yielding as corollary the "missing part" of our monograph (up to the mentioned exception). In particular we will classify the translation generalized quadrangles of order $(s, s^2)$ containing a regular line which is not incident with the translation point. As we start from a local situation the use of group theory is excluded. Hence, together with the results in the monograph of Payne and Thas \cite {P&JAT:84,P&JAT:09}, this will provide a geometric proof, relying just on a few elementary results on groups, of a large part of the deep group theoretical result of Fong and Seitz.
\par Also we mention that the classification of all Moufang generalized polygons, including the infinite ones, was obtained by Tits and Weiss, and this huge achievement is the theme of their monograph "Moufang Polygons" \cite {T&W:02}. In this context we also mention the excellent monograph on Moufang quadrangles by De Medts \cite {DM:05}.
\par Finally we mention that in the period 1994-1999 Thas wrote a series of three papers \cite {JAT:94,JAT:97,JAT:99} on generalized quadrangles of order $(s, s^2)$. This paper can be considered as Part IV of that series.

\section{Finite generalized quadrangles}

\subsection{Definitions} 
A (finite) {\em generalized quadrangle} (GQ) is an incidence  structure $\mathcal {S} = (P, B, \mathbf {I})$ in which $P$ and $B$ are disjoint (nonempty) sets of objects called {\em points} and {\em lines} (respectively) and for which $\mathbf {I}$ is a symmetric point-line incidence relation satisfying the following axioms:
\begin {itemize}
\item [(i)] Each point is incident with $1 + t$ lines $(t \geq 1)$ and two distinct points are incident with at most one line.
\item [(ii)] Each line is incident with $1 + s$ points $(s \geq 1)$ and two distinct lines are incident with at most one point.
\item [(iii)] If $x$ is a point and $L$ is a line not incident with $x$, then there is a unique pair $(y, M) \in P \times B$ for which $x \mathbf{I} M \mathbf{I} y \mathbf{I} L$.
\end {itemize}
\par The integers  $s$ and $t$ are the {\em parameters} of the GQ and $\mathcal {S}$ is said to have {\em order} $(s, t)$; if $s = t$, then $\mathcal{S}$ is said to have {\em order} $s$. There is a point-line duality for GQ for which in any definition or theorem the words "point" and "line" are interchanged and the parameters $s$ and $t$ are interchanged. Normally, we assume without further notice that the dual of a given theorem or definition has also been given.
\par Given two (not necessarily distinct) points $x, y$ of the GQ $\mathcal{S}$, we write $x \sim y$ and say that $x$ and $y$ are {\em collinear}, provided that there is some line $L$ for which $x \mathbf{I} L \mathbf{I} y$. And $x \not \sim y$ means that $x$ and $y$ are not collinear. Dually, for $L, M \in B$, we write $L \sim M$  or $L \not \sim M$ according as $L$ and $M$ are {\em concurrent} or nonconcurrent, respectively. The line (respectively point) which is incident with distinct collinear points $x, y$ (respectively distinct concurrent lines $L, M$) is denoted by $xy$ (respectively $LM$ or $L \cap M$).
\par For $x \in P$, put $x^\perp = \lbrace y \in P \ \vert\ y \sim x \rbrace$, and note that $x \in x^\perp$. If $A \subseteq P$, then A {\em "perp"} is defined by $A^\perp = \cap \lbrace x^\perp \ \vert \ x \in A \rbrace$. So for $x, y \in P$, with $x \not = y$, we have $\lbrace x, y \rbrace ^\perp = x^\perp \cap y^\perp$ and $\lbrace x, y \rbrace^{\perp \perp} = \lbrace u \in P \ \vert \ u\in z^\perp \ \forall z \in x^\perp \cap y^\perp \rbrace$. So $\vert \lbrace x, y \rbrace^\perp \vert = s + 1$ or $t + 1$ according as $x \sim y$ or $x \not \sim y$, and $\vert \lbrace x, y \rbrace ^{\perp \perp} \vert = s + 1$ or $\vert \lbrace x, y \rbrace ^{\perp \perp} \vert \leq t + 1$ according as $x \sim y$ or $x \not \sim y$.
\par A {\em triad} (of points) is a triple of pairwise noncollinear points. Given a triad $T = \lbrace x, y, z \rbrace$, a {\em center} of $T$ is just a point of $T^\perp$.
\par For terminology, notations, results, etc., concerning finite generalized quadrangles and not explicitly given here, see the monograph of Payne and Thas \cite {P&JAT:84,P&JAT:09}.

\subsection{Restrictions on parameters}
Let $\mathcal{S} = (P, B, \mathbf{I})$ be a GQ of order $(s, t)$. If $\vert P \vert = v$ and $\vert B \vert = b$, then $v = (s + 1)(st + 1)$ and $b = (t + 1)(st +1)$; see 1.2.1 of \cite {P&JAT:84,P&JAT:09}. Also, $s + t$ divides $st(s + 1)(t + 1)$; see 1.2.2 of \cite {P&JAT:84,P&JAT:09}.
\par If $s > 1$ and $t > 1$, then $t \leq s^2$, and dually $s \leq t^2$; see 1.2.3 of \cite {P&JAT:84,P&JAT:09}. Also, for $s > 1$ and $t > 1$, $s^2 = t$ if and only if each triad (of points) has a constant number of centers, in which case this constant number of centers is $s + 1$; see 1.2.4 of \cite {P&JAT:84,P&JAT:09}.

\subsection{Regularity}
Let $\mathcal{S} = (P,B, \mathbf{I})$ be a GQ of order $(s, t)$. If $x \sim y, x \not = y$, or if $x \not \sim y$ and $\vert \lbrace x, y \rbrace^{\perp \perp} \vert = t + 1$, we say that the pair $\lbrace x, y \rbrace$ is {\em regular}. The point $x$ is {\em regular} provided $\lbrace x, y \rbrace$ is regular for all $y \in P, y \not = x$. A point $x$ is {\em coregular} provided each line incident with $x$ is regular.
\par If $\mathcal{S}$ contains a regulair pair $\lbrace x, y \rbrace$, then either $s =1$ or $s \geq t$; see 1.3.6 of \cite {P&JAT:84,P&JAT:09}. A coregular point in a GQ of order $s, s > 1,$ is regular if and only $s$ is even; see 1.5.2 of \cite {P&JAT:84,P&JAT:09}.

\section{Translation generalized quadrangles}
\subsection{Elations, translations, symmetries}
Let $\mathcal{S} = (P, B, \mathbf{I})$ be a GQ of order $(s, t)$, with $s \not = 1, t \not = 1$. A collineation $\theta$ of $\mathcal{S}$ is an {\em elation} about the point $p$ if $\theta$ = id or if $\theta$ fixes all lines incident with $p$ and fixes no point of $P \setminus p^\perp$. If there is a group $H$ of elations about $p$ acting regularly on $P \setminus p^\perp$, we say $\mathcal{S}$ is an {\em elation generalized quadrangle} (EGQ) with {\em elation group} $H$ and {\em base point} $p$. Briefly we say that $(\mathcal{S}^{(p)}, H)$ or $\mathcal{S}^{(p)}$ is an EGQ. If the group $H$ is abelian, then we say that the EGQ $(\mathcal{S}^{(p)}, H)$ is a {\em translation generalized quadrangle} (TGQ) with {\em base point} or {\em translation point} $p$ and {\em translation group} $H$. In such a case $H$ is the set of all elations about $p$; see 8.6.4 of \cite {P&JAT:84,09}. For any TGQ $\mathcal{S}^{(p)}$ the point $p$ is coregular, hence the parameters $s$ and $t$ of a TGQ satisfy $s \leq t$; see Section 2.3. Also, by 8.5.2 of \cite {P&JAT:84,P&JAT:09}, for any TGQ with $s \not = t$ we have $s = q^a$ and $t  = q^{a + 1}$ with $q$ a prime power and $a$ an odd integer; if $s$ (or $t$) is even, then by 8.6.1(iv) of \cite {P&JAT:84,P&JAT:09} either $s = t$ or $s^2 = t$.
\par If a collineation $\theta$ of the GQ $\mathcal{S} = (P, B, \mathbf{I})$ of order $(s, t), s \not = 1 \not = t$, fixes each line of $L^\perp$, with $L \in B$, then $\theta$ is a {\em symmetry} about $L$. By 8.1 of \cite {P&JAT:84,09} any symmetry about $L$ is an elation about $L$. The line $L$ is called an {\em axis of symmetry} provided its group of symmetries has order $s$. If $\mathcal{S}$ is a TGQ with translation point $p$, then the lines $L_0, L_1, \ldots, L_t$ incident with $p$ are axes of symmetry and the translation group $G$ contains the $t + 1$ full groups of symmetries about these lines; in fact $G$ is generated by these symmetries; see 8.3.1 of \cite {P&JAT:84,09}. Conversely, if $p$ is a point of the GQ $\mathcal{S}$ for which $t + 1$ lines incident with $p$ are axes of symmetry, then $(\mathcal{S}^{(P}, G)$ is a TGQ, where $G$ is generated by the symmetries about the lines through $p$; see 8.3.1 of \cite {P&JAT:84,P&JAT:09}.
\par For a detailed study of TGQ we refer to the monograph "Translation Generalized Quadrangles" by Thas, Thas and Van Maldeghem \cite {TTVM:06}.

\subsection{The projective model}
In $\PG(2n + m - 1, q)$ consider a set $O = O(n, m, q)$ of $q^m + 1$ $(n - 1)$-dimensional subspaces $\PG^{(0)}(n - 1, q),
\PG^{(1)}(n - 1, q), \ldots, \PG^{(q^m)}(n - 1, q)$, every three of which generate a $\PG(3n - 1, q)$ and such that each element $\PG^{(i)}(n - 1, q)$ of $O(n, m, q)$ is contained in a $\PG^{(i)}(n + m - 1, q)$ having no point in common with any $\PG^{(j)}(n - 1, q)$, for $j \not = i$. It is easy to check that $\PG^{(i)}(n + m - 1, q)$ is uniquely determined, $i = 0, 1, \ldots, q^m$. The space $\PG^{(i)}(n + m - 1, q)$ is called the {\em tangent space} of $O(n, m, q)$ at $\PG^{(i)}(n - 1, q)$. For $n = m$ such a set $O(n, n, q)$ is called a {\em pseudo-oval} or a {\em generalized oval} or an {\em $\lbrack n - 1 \rbrack$-oval} of $\PG(3n - 1, q)$; a generalized oval of $\PG(2, q)$ is just an oval of $\PG(2, q)$. For $n \not = m$ such a set $O(n, m, q)$ is called a {\em pseudo-ovoid} or a {\em generalized ovoid} or an {\em $\lbrack n - 1 \rbrack$-ovoid} or an {\em egg} of $\PG(2n + m - 1, q)$; a $\lbrack 0 \rbrack$-ovoid of $\PG(3, q)$ is just an ovoid of $\PG(3, q)$.
\par Now embed $\PG(2n + m - 1, q)$ in a $\PG(2n + m, q)$ and construct a point-line geometry $T(n, m, q) = T(O)$ as follows. \\
Points are of three types:
\begin{itemize}
\item[(i)] the points of $\PG(2n + m, q) \setminus \PG(2n + m - 1, q)$;
\item[(ii)] the $(n + m)$-dimensional subspaces of $\PG(2n + m,q)$ wich intersect $\PG(2n + m - 1, q)$ in one of the spaces $\PG^{(i)}(n + m - 1, q)$;
\item[(iii)] the symbol $(\infty)$.
\end{itemize}
Lines are of two types:
\begin{itemize}
\item[(a)] the $n$-dimensional subspaces of $\PG(2n + m, q)$ which intersect $\PG(2n + m - 1, q)$ in one of the spaces $\PG^{(i)}(n - 1, q)$;
\item[(b)] the elements of $O(n, m, q)$.
\end{itemize}
\par Incidence in $T(n, m, q)$ is defined as follows. A point of type (i) is incident only with lines of type (a); here the incidence is that of $\PG(2n + m, q)$. A point of type (ii) is incident with all lines of type (a) contained in it and with the unique element of $O(n, m, q)$ contained in it. The point $(\infty)$ is incident with no line of type (a) and with all lines of type (b).

\begin{theorem} [8.7.1 of Payne and Thas \cite {P&JAT:84,P&JAT:09}]
\label{thm3.1}
The point-line geometry $T(O) = T(n, m, q)$ is a TGQ of order $(q^n, q^m)$ with base point $(\infty)$. Conversely, every TGQ is isomorphic to a $T(n, m, q)$. It follows that the theory of TGQ is equivalent to the theory of the sets $O(n, m, q)$.
\end{theorem}

\begin{corollary}
\label{cor3.2}
The following hold for every $O(n, m, q)$:
\begin{itemize}
\item[(i)] $n = m$ or $n(a + 1) = ma$ with $a$ odd;
\item[(ii)] if $q$ is even, then $n = m$ or $m = 2n$.
\end{itemize}
\end{corollary}

\par Each TGQ $\mathcal{S}$ of order $(s, t)$ with base point $(\infty)$, with $s \not =1 \not =t$, has a {\em kernel} $\mathbb {K}$,
where $\mathbb {K}$ is a field the multiplicative group of which is isomorphic to the group of all collineations of $\mathcal {S}$ fixing linewise the point $(\infty)$ and any given point not collinear with $(\infty)$; see 8.6.5 of Payne and Thas \cite {P&JAT:84,P&JAT:09}.
\par We have $\vert \mathbb {K} \vert \leq s$; see Payne and Thas \cite {P&JAT:84,P&JAT:09}. The field GF$(q)$ is a subfield of $\mathbb {K}$ if and only if $\mathcal {S}$ is isomorphic to a $T(n, m, q)$; see 8.7.1 of Payne and Thas \cite {P&JAT:84,P&JAT:09}.
\par Consider any GQ $T(n, m, q) = T(O)$, with $n = m$ and $q$ odd or with $n \not = m$. By 8.7.2 of \cite {P&JAT:84,P&JAT:09} or 3.9 and 3.10 of \cite {TTVM:06} the $q^m + 1$ tangent spaces of $O = O(n, m, q)$ form an $O^\star = O^\star(n, m, q)$ in the dual space of $\PG(2n + m - 1, q)$. So in addition to $T(O)$ there arises a TGQ $T(O^\star)$, also denoted $T^\star(O)$, with the same parameters. The TGQ $T^\star(O)$ is called the {\em translation dual} of the TGQ $T(O)$. Examples are known for which $T(O) \cong T^\star(O)$ and examples are known for which $T(O) \not \cong T^\star(O)$; see Thas \cite {JAT:94} and Thas, Thas and Van Maldeghem \cite {TTVM:06}.
\par Let $O(n, 2n, q)$ be an egg of $\PG(4n - 1, q)$. We say that $O(n, 2n, q)$ is {\em good} at the element $\PG^{(i)}(n - 1, q)$ of $O(n, 2n, q)$ if any $\PG(3n - 1, q)$ containing $\PG^{(i)}(n - 1, q)$ and at least two other elements of $O(n, 2n, q)$, contains exactly $q^n + 1$ elements of $O(n, 2n, q)$.

\section{Flocks and generalized quadrangles}
Let $F$ be a {\em flock} of the quadratic cone $K$ of $\PG(3, q)$ with vertex $x$, that is, a partition of $K\setminus \lbrace x \rbrace$ into $q$ disjoint irreducible conics. Independently, Walker \cite {W:76} and Thas \cite {JAT:76} show that with each flock there corresponds a translation plane $\mathcal{P}$ of order $q^2$. If $\mathcal{P}$ is a semifield plane, then $F$ is called a {\em semifield flock}; see e.g. \cite {TTVM:06}.
\par Relying on the work of Payne \cite {P:80,P:85} and Kantor \cite {K:86}, Thas \cite {JAT:87} proves that with $F$ there corresponds a GQ $\mathcal{S}(F)$ of order $(q^2, q)$. Such a GQ $\mathcal{S}(F)$ will be called a {\em flock GQ}.
\par For the following theorem one has to rely on work of Casse, Thas and Wild {\cite {CTW:85}, Lunardon \cite {L:97} and Payne \cite {P:89}; see also 4.7 and 4.8 of Thas, Thas and Van Maldeghem \cite {TTVM:06}.

\begin{theorem}
\label{thm4.1}
The flock GQ $\mathcal{S}(F)$ is a TGQ if and only if $F$ is a semifield flock.
\end{theorem}

\begin{theorem} [Thas \cite {JAT:94,JAT:99}]
\label{thm4.2}
If the egg $O$ is good, then the point-line dual of the translation dual $T(O^\star)$ of $T(O)$ is a flock GQ $\mathcal{S}(F)$.
\end {theorem}

\begin{corollary}
\label{cor4.3}
If the egg $O$ is good, then the point-line dual of the translation dual $T(O^\star)$ of $T(O)$ is a semifield flock GQ $\mathcal{S}(F)$.
\end{corollary}

\section{Moufang quadrangles}
\subsection{Moufang conditions}
Let $\mathcal{S} = (P, B, \mathbf{I})$ be a GQ of order $(s, t)$. For a fixed point $p$ define the following condition.

$(M)_p$: {\em For any two distinct lines $A$ and $B$ of $\mathcal{S}$ incident with $p$, the group of collineations of $\mathcal{S}$ fixing $A$ and $B$ pointwise and $p$ linewise is transitive on the set of lines distinct from $A$ and incident with a given point $x$ on $A$, with $x \not = p$}.

The GQ $\mathcal{S}$  is said to satisfy condition $(M)$ provided it satsifies $(M)_p$ for all points $p \in P$. For a fixed line $L \in B$ let $(\widehat M)_L$ be the condition that is the dual of $(M)_p$, and let $(\widehat M)$ be the dual of $(M)$. If $s \not = 1 \not = t$ and $\mathcal{S}$ satisfies both $(M)$ and $(\widehat M)$ it is said to be a {\em Moufang} GQ; if $s \not = 1 \not = t$ and $\mathcal{S}$ satisfies either $(M)$ or $(\widehat M)$ it is said to be a {\em half Moufang} GQ.
\par A celebrated result of Tits \cite {T:76} is that all finite Moufang GQ are classical or dual classical. In fact, Tits observes that his result easily follows from the classification of all groups with a $BN$-pair of rank 2 by Fong and Seitz \cite {F&S:73,F&S:74}. A huge achievement of Tits and Weiss \cite {T&W:02} is the classification of all (finite and infinite) Moufang GQ. In \cite {TPVM:91} Thas, Payne and Van Maldeghem prove that half Moufang implies Moufang for finite GQ; in \cite {Te:04} Tent proves the analogue for infinite GQ. 
\par In their monograph Payne and Thas \cite {P&JAT:84,P&JAT:09} did an almost successful attempt to prove the Moufang theorem (in the finite case) in a pure geometric way. To complete a geometric proof of the theorem of Tits (in the finite case) it would be sufficient to show in a geometric way that if, for each point $p$ of $\mathcal{S}$, the GQ $\mathcal{S}^{(p)}$ is a TGQ of order $(s, s^2), s \not = 1$, then $\mathcal{S}$ is isomorphic to the classical GQ $Q(5, s)$ arising from a non-singular elliptic quadric in $\PG(5, s)$. As a corollary there would also follow a pure geometric proof of a large part of the theorem of Fong and Seitz, namely in the case where the Weyl group is $D_8$.
\par That "missing part" will be a corollary of the much stronger results that follow.

\section{Translation generalized quadrangles of order $(s, s^2)$ with a regular line not containing the base point}
\subsection{Introduction}
In this section let $\mathcal{S} = (P, B, \mathbf{I})$ be a TGQ of order $(s, s^2), s \not = 1$, with base point $(\infty)$, and having a regular line $M$ not incident with $(\infty)$. Assume that $\mathcal{S} = T(n, 2n, q) = T(O)$, with $O = \lbrace \pi_0, \pi_1, \ldots, \pi_{q^{2n}} \rbrace$, and where $M$ is concurrent with the line $\pi_0$ of $\mathcal{S}$. By transitivity it is clear that all lines concurrent with $\pi_0$ are regular. The projective space containing $O$ will be denoted by $\PG(4n - 1, q)$ and the projective space containing $T(O)$ by $\PG(4n, q)$. The tangent space of $O$ at $\pi_i$ will be denoted by $\tau_i$, with $i = 0, 1, \ldots, q^{2n}$. Further, let $\widetilde O = \pi_0 \cup \pi_1 \cup \ldots \cup \pi_{q^{2n}}$.

\subsection{Conics on $\widetilde O$}
Consider lines $L_0, L_1$ of $\mathcal{S}$, with $L_0 \! \not \! \mathbf{I}  (\infty)  \! \not \! \mathbf{I} L_1, L_0 \sim \pi_0, L_1 \sim \pi_1,  L_0 \not \sim L_1$. Let $\lbrace L_0, L_1 \rbrace ^\perp = \lbrace M_0, M_1, \ldots, M_{q^{n}} \rbrace, \lbrace L_0, L_1 \rbrace^{\perp \perp} = \lbrace L_0, L_1, \ldots, L_{q^{n}} \rbrace$ and assume that notations are chosen in such a way that $L_i \sim \pi_i \sim M_i$, with $i = 0, 1, \ldots q^n$. Then $\pi_i, L_i, M_i$ contain a common point $r_i$ of $\mathcal{S}$, with $r_i$ a $3n$-dimensional space containing $\tau_i, i = 0, 1, \ldots, q^n$. So the $(n + 1)$-dimensional space $\langle L_i, M_i \rangle$ generated by $L_i$ and $M_i$ intersects $\PG(4n - 1,q )$ in an $n$-dimensional space $\gamma_i$ containing $\pi_i$ and contained in $\tau_i$.
\par As $\langle L_i, M_i \rangle \cap \langle L_j, M_j \rangle$, with $i \not = j$ and $i, j \in \lbrace 0, 1, \ldots, q^n \rbrace$, is a line $U_{ij} = U_{ji}$ we have $\langle L_i, M_i \rangle \cap \langle L_j, M_j \rangle \cap \PG(4n - 1, q) = \gamma_i \cap \gamma_j = U_{ij} \cap \PG(4n - 1, q) = \lbrace u_{ij} \rbrace$. Hence $\gamma_i$ and $\gamma_j$ have exactly one point $u_{ij} = u_{ji}$ in common. This point $u_{ij}$ belongs to $\tau_i \cap \tau_j$.
\par The grid, that is, the subquadrangle of $\mathcal{S}$ of order $(s, 1)$, containing $L_1, L_0$ will be denoted by $\mathcal{G}$. Let us consider the lines $L_0, L_1, L_2, M_0, M_1, M_2$ of $\mathcal{G}$, and put $L_i \cap M_j = n_{ij}$, with $i \not = j$ and $ i, j \in \lbrace 0, 1, 2 \rbrace$. The space $\langle L_0, L_1 \rangle \cap \PG(4n - 1, q) = \langle M_0, M_1 \rangle \cap \PG(4n - 1, q) = \rho$ is the $2n$-dimensional space $\langle \gamma_0, \gamma_1 \rangle$ containing $\langle \pi_0, \pi_1 \rangle$. The line $n_{20}n_{21}$ of $\PG(4n, q)$ is contained in  $\langle L_0, L_1 \rangle$, so intersects $\rho$ is a point $s_2$. Notice that $\langle L_0, L_1 \rangle \cap L_2 = n_{20}n_{21}$. Also $\langle L_0, L_1, L_2, M_0, M_1, M_2 \rangle$ is $3n$-dimensional. The point $s_2$ is the unique point of $\rho$ in $\pi_2$. Hence also the line $n_{02}n_{12}$ contains $s_2$, and so $n_{20}n_{21} \cap n_{02}n_{12} = \lbrace s_2 \rbrace$. Analogously, let $n_{01}n_{02} \cap n_{10}n_{20} = \lbrace s_0 \rbrace$ and $n_{10}n_{12} \cap n_{01}n_{21} = \lbrace s_1 \rbrace$. Consequently the lines $n_{01}n_{02}, n_{10}n_{12}, n_{20}n_{21}, n_{10}n_{20}, n_{01}n_{21}, n_{02}n_{12}$ belong to a common hyperbolic quadric $\mathcal H$ of some $\PG(3, q)$.
\par Let $n \in n_{20}n_{21} \setminus \lbrace n_{20}, n_{21}, s_2 \rbrace$. Through $n$ there is a line $M$ of the GQ intersecting $L_0$ and $L_1$. Let $m_i$ be the common point of $M$ and $L_i$, $i = 0, 1$. Assume, by way of contradiction, that $n, m_0, m_1$ are not collinear in $\PG(4n, q)$. Then $M \cap \rho$ contains the line $\langle n, m_0, m_1 \rangle \cap \PG(4n - 1, q)$ (as $n, m_0, m_1$ belong to $\langle L_0, L_1 \rangle = \langle M_0, M_1 \rangle$), so $M$ contains a point of $\langle \pi_0, \pi_1 \rangle$, so $M \cap \PG(4n - 1, q) = \pi_i, i \not = 0, 1$, contains a point of $\langle \pi_0, \pi_1 \rangle$, clearly a contradiction. Hence the points $n, m_0, m_1$ are collinear in $\PG(4n, q)$. Let $n_3, n_4, \ldots, n_q$ be the points of $n_{20}n_{21}$ distinct from $n_{20}, n_{21}, s_2$, let $m_0^i$ be the point of $L_0$ collinear in $\mathcal{S}$ with $n_i$, let $m_1^i$ be the point of $L_1$ collinear in $\mathcal{S}$ with $n_i$, and let $m_0^im_1^i \cap \rho = \lbrace s_i \rbrace$, with $i = 3, 4, \ldots, q$. Notations are chosen in such a way that $s_i \in \pi_i$, with $i = 3, 4, \ldots, q$. From the foregoing it follows that the $q + 1$ lines of $\mathcal{H}$ which intersect $L_0, L_1, L_2$ are contained in the $n$-dimensional spaces $M_0, M_1, \ldots, M_q$ which respectively contain $\pi_0, \pi_1, \ldots, \pi_q$. Similarly the $q + 1$ lines of $\mathcal{H}$ intersecting $M_0, M_1, M_2$ are contained in respectively $L_0, L_1, \ldots, L_q$. Let $N_i$ be the line of $\mathcal{H}$ contained in $L_i$, and let $R_i$ be the line of $\mathcal{H}$ contained in $M_i, i= 0, 1, \ldots, q$. Then $N_i \cap R_i = \lbrace s_i \rbrace, i = 0, 1, \ldots, q$. So $\lbrace s_0, s_1, \dots, s_q \rbrace$ is a non-singular conic $\mathcal{C}$ contained in $\widetilde{O}$.
\par The tangent line $T_i$ of $\mathcal{C}$ at $s_i$ is contained in $\langle L_i, M_i \rangle \cap \PG(4n - 1, q)$, it is the line $\langle R_i, N_i \rangle \cap \PG(4n - 1, q)$. If $\sigma$ is the plane containing $\mathcal{C}$, then $T_i = \gamma_i \cap \sigma = \tau_i \cap \sigma$, with $i = 0, 1, \ldots, q$.
\par In the foregoing the lines $L_0$ and $L_1$ may be replaced by any two distinct lines of $\lbrace L_0, L_1 \rbrace ^{\perp \perp}$.

\begin{lemma}
\label{lem6.1}
The conic $\mathcal{C}$ is uniquely defined by $\pi_0, \pi_1$ and $s_2$.
\end{lemma}
{\em Proof}\quad
We have $T_i = \gamma_i \cap \sigma = \tau_i \cap \sigma$, with $\gamma_i = \langle L_i, M_i \rangle \cap \PG(4n - 1, q), i = 0, 1, \ldots, q$. Also, $\rho = \langle \gamma_0, \gamma_1 \rangle, \rho = \langle L_0, L_1 \rangle \cap \PG(4n - 1, q) = \langle \pi_0, \pi_1, s_2 \rangle, \gamma_i = \rho \cap \tau_i, i = 0, 1$. So $\gamma_0, \gamma_1$ and also $T_0 \cap T_1 = \gamma_0 \cap \gamma_1 = \lbrace u_{01} \rbrace$ are uniquely defined by $\pi_0, \pi_1$ and $s_2$.
\par Now we project from $u_{01}$ onto a $(2n - 1)$-dimensional space $\bar \rho$ of $\rho$, not containing $u_{01}$. Let the image of $\pi_i$ be $\bar \pi_i$, with $i = 0, 1$, and let the image of $s_2$ be $\bar s_2$. There is a unique line $\bar U$ in $\bar \rho$ containing $\bar s_2$ and intersecting $\bar \pi_0$ and $\bar \pi_1$. Then  $\langle u_{01}, \bar U \rangle$ is the unique plane in $\rho$ containing $u_{01}, s_2$ and intersecting $\pi_0$ and $\pi_1$. Hence $\langle u_{01}, \bar U \rangle = \sigma$ and $\lbrace s_i \rbrace = \sigma \cap \pi_i, i = 0, 1$. Consequently $s_0$ and $s_1$ are uniquely defined by $\pi_0, \pi_1$ and $s_2$. So $\mathcal{C}$ is the unique conic containing $s_0, s_1, s_2$ and having $s_0u_{01}$ and $s_1u_{01}$ as tangent lines. \eop \\

\begin{remark}
\label{rem6.2}
In \ref{lem6.1} the subscripts $0, 1, 2$ may be replaced by any three distinct subscripts $i, j, k \in \lbrace 0, 1, \ldots, q^n \rbrace$.
\end{remark}

{\bf Definition}.
A conic $\mathcal{C}$ of the type described in \ref{lem6.1} or in \ref{rem6.2}, where we start from any three distinct lines of some $\lbrace L_0, L_1 \rbrace^{\perp \perp}$, will be called a {\em $\widetilde \pi_0$-conic} of $\widetilde{O}$. If $\mathcal{C}$ intersects $\pi_0$ it will also be called a {\em $\pi_0$-conic} of $\widetilde{O}$.

\begin{lemma}
\label{lem6.3}
Let $\pi_i, \pi_j, \pi_k$ be distinct elements of $O$, with $\pi_0 \in \lbrace \pi_i, \pi_j, \pi_k \rbrace$, and let $s$ be a point of $\pi_k$. Then $s$ is contained in exactly one $\pi_0$-conic intersecting $\pi_i$ and $\pi_j$.
\end{lemma}
{\em Proof}\quad
Let $L_i$ and $L_j$ be nonconcurrent lines of $\mathcal{S}$ which are concurrent with respectively $\pi_i$ and $\pi_j$, but not incident with the point $(\infty)$. Assume also that $L_i$ and $L_j$ are chosen in such a way that $\langle L_i, L_j \rangle \cap \PG(4n - 1, q) = \langle \pi_i, \pi_j, s \rangle$. If $\pi_0 \in \lbrace \pi_i, \pi_j \rbrace$, then the lines $L_i$ and $L_j$ define a $\pi_0$-conic containing $s$ and intersecting $\pi_i$ and $\pi_j$. By \ref{lem6.1} this $\pi_0$-conic is unique. Now assume that $\pi_0 = \pi_k$. As in $\PG(4n, q)$ there is a line containing $s$ and intersecting $L_i$ and $L_j$, there is also a line $M$ in $\mathcal{S}$ concurrent with $\pi_k$ and intersecting $L_i$ and $L_j$. As the line $M$ is regular in $\mathcal{S}$, also the pair $\lbrace L_i, L_j \rbrace$ is regular in $\mathcal{S}$. Interchanging roles of $L_0$ and $L_i$, see also \ref{rem6.2}, the result follows. \eop \\

\subsection{Segre varieties defined by $\lbrace O, \pi_0 \rbrace$}
For Segre varieties we refer to the monograph "General Galois Geometries" by Hirschfeld and Thas \cite {H&JAT:91,H&JAT:16}.

\begin{lemma}
\label{lem6.4}
If $\pi_i$ and $\pi_j$ are distinct elements of $O \setminus \lbrace \pi_0 \rbrace$, then all $\pi_0$-conics intersecting $\pi_0, \pi_i, \pi_j$ belong to a common Segre variety $\mathcal{S}_{2; n - 1}$.
\end{lemma} 
{\em Proof}\quad
Let $\langle \pi_0, \pi_i, \pi_j \rangle = \PG(3n - 1, q)$ and $\tau_k \cap \PG(3n - 1, q) = \bar \tau_k$, with $k \in \lbrace 0, i, j \rbrace$. Then $\bar \tau_0, \bar \tau_i, \bar \tau_j$ are $(2n - 1)$-dimensional. Further, let $\bar \tau_i \cap \bar \tau_j = \eta$. So $\eta$ is $(n - 1)$-dimensional. Let $s$ be any point of $\pi_0$, and let $\mathcal{C}$ be the $\pi_0$-conic containing $s$ and intersecting $\pi_i$ and $\pi_j$. Let $\mathcal{C} \cap \pi_i = \lbrace s_i \rbrace$ and $\mathcal{C} \cap \pi_j = \lbrace s_j \rbrace$. The tangent lines of $\mathcal{C}$ at respectively $s_i$ and $s_j$ have a point $u_{ij}$ in common, which is the common point of $\langle s, \pi_i, \pi_j \rangle \cap \tau_i = \gamma_i$ and $\langle s, \pi_i, \pi_j \rangle \cap \tau_j = \gamma_j$; so $u_{ij}$ is contained in $\bar \tau_i \cap \bar \tau_j = \eta$. Hence the plane $\sigma$ containing $\mathcal{C}$ intersects the spaces $\pi_0, \pi_i, \pi_j, \eta$.
\par It follows that the $(q^n - 1)/(q - 1)$ planes $\sigma$ are exactly the $(q^n - 1)/(q - 1)$ planes of $\PG(3n - 1, q)$ which intersect $\pi_0, \pi_i, \pi_j, \eta$. It follows that the union of these planes is a Segre variety $\mathcal{S}_{2; n - 1}$ of a plane and an $(n - 1)$-dimensional space; see \cite {H&JAT:91,H&JAT:16}.
\par We conclude that the $(q^n - 1)/(q - 1)$ $\pi_0$-conics intersecting $\pi_0, \pi_i$ and $\pi_j$ belong to $\mathcal{S}_{2; n - 1}$. \eop \\

\begin{lemma}
\label{lem6.5}
Let $\pi_i$ and $\pi_j$ be distinct elements of $O \setminus \lbrace \pi_0 \rbrace$. Further, let $\langle \pi_0, \pi_i, \pi_j \rangle = \PG(3n - 1, q)$, and $\tau_l \cap \PG(3n - 1, q) = \bar \tau_l$, with $l \in \lbrace 0, i, j \rbrace$. Finally, let $\bar \tau_0 \cap \bar \tau_j = \eta_i, \bar \tau_0 \cap \bar \tau_i = \eta_j, \bar \tau_i \cap \bar \tau_j = \eta_0$. If $q$ is even, then $\eta_0 = \eta_i = \eta_j \subset \PG(3n - 1, q)$; if $q$ is odd, then $\langle \pi_0, \eta_0 \rangle, \langle \pi_i, \eta_i \rangle, \langle \pi_j , \eta_j \rangle$ have an $(n - 1)$-dimensional space in common.
\end{lemma}
{\em Proof}\quad
Let $q$ be even. By the proof of \ref{lem6.4} the $(q^n - 1)/(q - 1)$  nuclei of the corresponding $\pi_0$-conics intersecting $\pi_0, \pi_i, \pi_j$ belong to $\eta_0$, respectively $\eta_i, \eta_j$. Hence $\eta_0 = \eta_i = \eta_j$. Clearly that space belongs to $\PG(3n - 1, q)$.
\par Next, let $q$ be odd. Let $\mathcal{C}$ be a $\pi_0$-conic intersecting $\pi_0, \pi_i, \pi_j$, let $\mathcal{C} \cap \pi_l = \lbrace s_l \rbrace$ and let $T_l$ be the tangent line of $\mathcal{C}$ at $s_l$, with $l = 0, i, j$. Further, let $T_0 \cap T_i = \lbrace u_{0i} \rbrace, T_0 \cap T_j = \lbrace u_{0j} \rbrace, T_i \cap T_j = \lbrace u_{ij} \rbrace$. Then, by an elementary property of conics, the lines $s_iu_{0j}, s_ju_{oi}, s_0u_{ij}$ have a point $z$ in common. As the $(q^n - 1)/(q - 1)$ points $s_0$, respectively $s_i, s_j, u_{0i}, u_{0j}, u_{ij}$ belong to a maximal $(n - 1)$-space of $\mathcal{S}_{2; n - 1}$, also the $(q^n - 1)/(q - 1)$ points $z$ will belong to a maximal $(n - 1)$-space of $\mathcal{S}_{2; n - 1}$. It follows that the $(2n - 1)$-dimensional spaces $\langle \pi_0, \eta_0 \rangle, \langle \pi_i, \eta_i \rangle, \langle \pi_j, \eta_j \rangle$ have an $(n - 1)$-dimensional space in common.  \eop \\

\par Consider two distinct elements of $O \setminus \lbrace \pi_0 \rbrace$, say $\pi_1$ and $\pi_2$. Let $\langle \pi_0, \pi_1, \pi_2 \rangle = \PG(3n - 1, q	)$. The $\pi_0$-conics intersecting $\pi_0, \pi_1, \pi_2$ contain $(q^n - 1)(q + 1)/(q - 1)$ points of $\widetilde O$. By \ref{lem6.4} these points belong to a common Segre variety $\mathcal{S}_{2; n - 1}$. Let $\mathcal{C}$ be some $\pi_0$-conic intersecting $\pi_0, \pi_1, \pi_2$ and let $\mathcal{V}$ be the set of all points on the $(q^n - 1)/(q - 1)$ $\pi_0$-conics intersecting $\pi_0, \pi_1, \pi_2$. Further, let $s \in \mathcal{C}$.

\begin{lemma}
\label{lem6.6}
The maximal $(n - 1)$-space of $\mathcal{S}_{2; n - 1}$ (distinct from the plane of $\mathcal{C}$ if $n = 3$) containing $s$ is contained in $\mathcal{V}$.
\end{lemma}
{\em Proof}\quad
Let $\sigma$ be the plane of $\mathcal{C}$ and let $\sigma^\prime$ be the plane of a second $\pi_0$-conic $\mathcal{C}^\prime$ intersecting $\pi_0, \pi_1, \pi_2$. The maximal $(n - 1)$-spaces of $\mathcal{S}_{2; n - 1}$ (of the family containing $\pi_0, \pi_1, \pi_2$ if $n = 3$) define a linear projectivity $\theta$ from $\sigma$ onto $\sigma^\prime$. If $\mathcal{C} \cap \pi_i = \lbrace s_i \rbrace, \mathcal{C}^\prime \cap \pi_i  = \lbrace s^\prime_i \rbrace$, then $s_i^\theta = s^\prime_i, i = 0, 1, 2$. Let $\bar \tau_i = \tau_i \cap \langle \pi_0, \pi_1, \pi_2 \rangle$, with $i = 0, 1, 2$. If $\bar \tau_i \cap \bar \tau_j = \eta_k$ with $\lbrace i, j, k \rbrace = \lbrace 0, 1, 2 \rbrace$, and $\sigma \cap \eta_k = \lbrace l_k \rbrace, \sigma^\prime \cap \eta_k = \lbrace l^\prime_k \rbrace, k = 0, 1, 2$, then $l_k$ is the intersection of the tangent lines of $\mathcal{C}$ at $s_i$ and $s_j$, and $l^\prime_k$ is the intersection of the tangent lines of $\mathcal{C}^\prime$ at $s^\prime_i$ and $s^\prime_j$, with $\lbrace i, j, k \rbrace = \lbrace 0, 1, 2 \rbrace$. Also $\l^\theta_k = l^\prime_k$, with $k = 0, 1, 2$. It follows that $\mathcal{C}^\theta = \mathcal{C}^\prime$. Hence for any point $s \in \mathcal{C}$ the maximal $(n - 1)$-space of $\mathcal{S}_{2; n - 1}$ (distinct from the plane of $\mathcal{C}$ if $n = 3$) containing $s$ is contained in $\mathcal{V}$. \eop\\

\begin{lemma}
\label{lem6.7}
The set $\mathcal{V}$ is the union of $q + 1$ elements of $O$.
\end{lemma}
{\em Proof} \quad
From \ref{lem6.6} follows that $\mathcal{V}$ is the union of $q + 1$ $(n - 1)$-dimensional spaces. Let $\pi$ be such an $(n - 1)$-dimensional space. As $\pi \subset \widetilde{O}$, any line of $\pi$ belongs to some element of $O$. Hence $\pi \in O$, so $\mathcal{V}$ is the union of $q + 1$ elements of $O$. \eop \\

{\bf Definition}. The set $C$ consisting of the $q + 1$ elements of $O$ which intersect some $\pi_0$-conic $\mathcal{C}$ of $\widetilde O$ is called a {\em $\Pi_0$-conic} of $O$. Also, $\langle C \rangle \cap \tau_i = \bar \tau_i$, with $\pi_i \in C$, is called the {\em tangent space} of $C$ at $\pi_i$. Clearly any three distinct elements of $O$ are contained in exactly one $\Pi_0$-conic.
\\

\par Let $\delta$ be the projection from $\pi_0$ onto a $(3n - 1)$-dimensional subspace $\Phi$ of $\PG(4n - 1, q)$ which is skew to $\pi_0$. Let $\tau_0 \cap \Phi = \tau^\star_0$ and let $\langle \pi_0, \pi_j \rangle \cap \Phi = \pi^\star_j$, with $j = 1, 2, \ldots, q^{2n}$. Then $\tau^\star_0 \cup \pi^\star_1 \cup \pi^\star_2 \cup \ldots \cup \pi^\star_{q^{2n}} = \Phi$. Further, let $\langle \pi^\star_i, \pi^\star_j \rangle \cap \tau^\star_0 = \pi^\star_{ij} = \pi^\star_{ji}$, with $0 \not = i \not = j \not = 0$. The space $\tau^\star_0$ is $(2n - 1)$-dimensional, the spaces $\pi^\star_j$ are $(n - 1)$-dimensional, and the spaces $\pi^\star_{ij}$ are $(n - 1)$-dimensional.

\begin{lemma}
\label{lem6.8}
For any $i, j$, with $0 \not = i \not = j \not = 0$, the Segre variety $\mathcal{S}_{1; n - 1}$ defined by $\pi^\star_i, \pi^\star_j, \pi^\star_{ij}$, contains exactly $q$ elements $\pi^\star_k$.
\end{lemma}
{\em Proof} \quad
Consider a $\pi_0$-conic $\mathcal{C}$ intersecting $\pi_0, \pi_i, \pi_j$. The $q + 1$ elements of the $\Pi_0$-conic $C$ defined by $\mathcal{C}$ are maximal $(n - 1)$-spaces of a Segre variety $\mathcal{S}_{2; n - 1}$; see \ref{lem6.6} and \ref{lem6.7}. Let $\langle \pi_0, \pi_i, \pi_j \rangle \cap \tau_0 = \bar \tau_0$ and let $\pi_0, \pi_i, \pi_j, \pi_{k_1}, \ldots, \pi_{k_{q - 2}}$ be the elements of the $\Pi_0$-conic $C$. Then $\bar \tau_0 \cap \Phi = \pi^\star_{ij}, \pi^\star_i, \pi^\star_j, \pi^\star_{k_1}, \ldots, \pi^\star_{k_{q - 2}}$ are the elements of a system of maximal $(n - 1)$-dimensional spaces of a $\mathcal{S}_{1; n - 1}$. So the Segre variety $\mathcal{S}_{1; n - 1}$ defined by $\pi^\star_i, \pi^\star_j, \pi^\star_{ij}$ contains exactly $q$ elements $\pi^\star_k, k \in \lbrace 1, 2, \ldots, q^{2n}\rbrace$. \eop \\

\begin{remark}
\label{rem6.9}
For $n \not = 2$ the elements of the family of maximal $(n - 1)$-dimensional spaces of a $\mathcal{S}_{1; n - 1}$, and for $n = 2$ the elements of any of the two families of lines of $\mathcal{S}_{1; 1}$, form an {\em $(n - 1)$-regulus} of $\langle \mathcal{S}_{1; n - 1} \rangle$; see Chapter 4 of \cite {H&JAT:16}.
\end{remark}

\subsection{$\Pi_0$-sets $(\Pi, \Gamma)$}
Consider again $\lbrace L_0, L_1 \rbrace^{\perp \perp} = \lbrace L_0, L_1, \ldots, L_{q^n} \rbrace, L_0 \! \not \! \mathbf{I}  (\infty) \! \not \mathbf{I} L_1$, with $L_i \sim \pi_i, i = 0, 1, \ldots, q^n$. Also, let $\lbrace L_0, L_1 \rbrace^\perp = \lbrace M_0, M_1, \ldots, M_{q^n} \rbrace$, with $\pi_i \sim M_i \sim L_i \sim \pi_i, i = 0, 1, \ldots, q^n$. Further, $\langle L_i, M_i \rangle \cap \tau_i = \gamma_i, i = 0, 1, \ldots, q^n$. Then $\pi_i \subset \gamma_i$.
\\

{\bf Definition}. With $\lbrace \pi_0, \pi_1, \ldots, \pi_{q^n} \rbrace = \Pi$ and $\lbrace \gamma_0, \gamma_1, \ldots, \gamma_{q^n} \rbrace = \Gamma$, the pair $(\Pi, \Gamma)$ is called a {\em $\Pi_0$-set}.

\begin{lemma}
\label{lem6.10}
\begin{itemize}
\item [(i)] $\gamma_i \cap \gamma_j = \lbrace u_{ij} \rbrace$, with $i, j$ distinct elements of $\lbrace 0, 1, \ldots, q^n \rbrace$.
\item [(ii)] Any three distinct elements $\pi_i, \pi_j, \pi_k$, with $0 \in \lbrace i, j, k \rbrace \subset \lbrace 0, 1, \ldots, q^n \rbrace$, are contained in exactly one $\Pi_0$-conic of $\Pi$.
\item [(iii)] If $\pi_i, \pi_j, \pi_k$, with $0 \in \lbrace i, j, k \rbrace$, are distinct elements of $\Pi$, then $\vert \gamma_i \cap \langle \pi_j, \pi_k \rangle \vert = 1$.
\item [(iv)] If we know $\Pi$ and some $\gamma_i \in \Gamma$, then all $\gamma_j$ of $\Gamma$ are determined.
\end{itemize}
\end{lemma}
{\em Proof} \quad
\begin{itemize}
\item [(i)] Was proved in Section 6.2.
\item [(ii)] From \ref{lem6.3}, \ref{lem6.4}, \ref{lem6.6}, \ref{lem6.7} follows (ii).
\item [(iii)] The spaces $\pi_i, \pi_j, \pi_k$, with $0 \in \lbrace i, j, k \rbrace$, are contained in exactly one $\Pi_0$-conic $C$. The tangent space $\bar \tau_i$ of $C$ at $\pi_i$ contains the space $\gamma_i$. Hence the $n$-space $\gamma_i$ has exactly one point in common with the $(2n - 1)$-space $\langle \pi_j, \pi_k \rangle$ of the $(3n - 1)$-space $\langle C \rangle$.
\item [(iv)] Assume we know $\Pi$ and some $\gamma_i \in \Gamma$. Let $\pi_j \in \Pi, i \not = j$. Then $\tau_j$ has exactly one point $x$ in common with $\gamma_i$. So $\gamma_j = \langle \pi_j, x \rangle.$
\end{itemize} \eop \\

\begin{lemma}
\label{lem6.11}
The $\Pi_0$-set $(\Pi, \Gamma)$ is uniquely defined by $\lbrace \gamma_0, \pi_j \rbrace$ or $\lbrace \gamma_j, \pi_0 \rbrace$, with $j \in \lbrace 1, 2, \ldots, q^n \rbrace$. The elements of $\Pi \setminus \lbrace \pi_0, \pi_j \rbrace$ are the $q^n - 1$ elements of $O \setminus \lbrace \pi_0, \pi_j \rbrace$ containing the common points of $\widetilde O \setminus (\pi_0 \cup \pi_j)$ and $\langle \gamma_0, \pi_j \rangle$.
\end{lemma}
{\em Proof} \quad
By the proof of \ref{lem6.10}(iv) it is sufficient to assume that $\lbrace \gamma_0, \pi_j \rbrace$ is given. Each $n$-dimensional space containing $\pi_0$, distinct from $\gamma_0$, and contained in $\langle \gamma_0, \pi_j \rangle$ has exactly one point in common with $\widetilde O \setminus \pi_0$. In this way there arise $q^n - 1$ points of $\widetilde O \setminus (\pi_0 \cup \pi_j)$. Consider $\pi_k \in \Pi, \pi_k \not \in \lbrace \pi_0, \pi_j \rbrace$. Then by \ref{lem6.10} $\gamma_0 \subset \langle \pi_0, \pi_j, \pi_k \rangle$, so $\langle \gamma_0, \pi_j \rangle$ has exactly one point in common with $\pi_k$. Hence the $q^n - 1$ common points of $\widetilde O \setminus (\pi_0 \cup \pi_j)$ and $\langle \gamma_0, \pi_j \rangle$ determine the $q^n - 1$ elements of $\Pi \setminus \lbrace \pi_0, \pi_j \rbrace$. \eop\\

\begin{lemma}
\label{lem6.12}
Let $\gamma_0$ be any $n$-dimensional space containing $\pi_0$ and let $\pi_k \in O \setminus \lbrace \pi_0 \rbrace$. Then $\lbrace \gamma_0, \pi_k \rbrace$ determines exactly one $\Pi_0$-set $(\Pi, \Gamma)$.
\end{lemma}
{\em Proof} \quad
Let $L_0, M_0$ be two $n$-dimensional spaces of $\PG(4n, q)$ for which $\langle L_0, M_0 \rangle \cap \PG(4n - 1, q) = \gamma_0$. Let $L_k$ be an $n$-dimensional space containing $\pi_k$ and a point of $M_0 \setminus \pi_0$. Then the $\Pi_0$-set $(\Pi, \Gamma)$ corresponding to $\lbrace L_0, L_k \rbrace^{\perp \perp}$ is uniquely determined by the pair $\lbrace \gamma_0, \pi_k \rbrace$. \eop \\

\begin{lemma}
\label{lem6.13}
Let $\gamma_0$ be any $n$-dimensional space containing $\pi_0$. Then the $\Pi_0$-sets $(\Pi, \Gamma)$ with $\gamma_0 \in \Gamma$ determine a partition of $O \setminus \lbrace \pi_0 \rbrace$.
\end{lemma}
{\em Proof} \quad
By \ref{lem6.11} and \ref{lem6.12} all $\Pi_0$-sets $(\Pi, \Gamma)$ with $\gamma_0 \in \Gamma$ give rise to sets $\Pi \setminus \lbrace \pi_0 \rbrace$ defining a partition of $O \setminus \lbrace \pi_0 \rbrace$. \eop \\

\begin{remark}
\label{rem6.14}
Consider any $\Pi_0$-set $(\Pi, \Gamma)$. From Section 6.3 follows that for $q$ even all elements of $\Gamma$ have a point $n$ in common. This point $n$ is called the {\em kernel} of $(\Pi, \Gamma)$. The tangent spaces of $O$ containing the elements of $\Pi$ are exactly the $q^n + 1$ tangent spaces of $O$ containing the point $n$.
Now let $q$ be odd. Consider $\pi_0, \pi_i, \pi_j \in \Pi$, with $0 \not =  i \not = j \not = 0$. Let $\gamma_0 \cap \gamma_i = \lbrace u_{0i} \rbrace, \gamma_0 \cap \gamma_j = \lbrace u_{0j} \rbrace, \gamma_i \cap \gamma_j = \lbrace u_{ij} \rbrace$. Then the plane $\sigma = \langle u_{0i}, u_{0j}, u_{ij} \rangle$ contains a $\pi_0$-conic $\mathcal{C}$ and let $\mathcal{C} \cap \pi_0 = \lbrace s_0 \rbrace, \mathcal{C} \cap \pi_i = \lbrace s_i \rbrace, \mathcal{C} \cap \pi_j = \lbrace s_j \rbrace$. Then the lines $u_{0i}s_j, u_{0j}s_i$ and $u_{ij}s_0$ have a point $x$ in common.
\end{remark}

\subsection{Projections and projectivities}
In Section 6.3 we introduced the projection $\delta$ from $\pi_0$ onto a $(3n - 1)$-dimensional subspace $\Phi$ of $\PG(4n - 1, q)$ which is skew to $\pi_0$. Let $\tau_0 \cap \Phi = \tau^\star_0 = (\tau_0 \setminus \pi_0)^\delta$ and let $\langle \pi_0, \pi_j \rangle \cap \Phi = \pi^\star_j = \pi^\delta_j$, with $j = 1, 2, \ldots, q^{2n}$. Then $\tau^\star_0 \cup \pi^\star_1 \cup \pi^\star_2 \cup \ldots \cup \pi^\star_{q^{2n}} = \Phi$. Further, let $\langle \pi^\star_i, \pi^\star_j \rangle \cap \tau^\star_0 = \pi^\star_{ij} = \pi^\star_{ji}$, with $0 \not = i \not = j \not = 0$. The space $\tau^\star_0$ is $(2n - 1)$-dimensional, the spaces $\pi^\star_j$ are $(n - 1)$-dimensional, and the spaces $\pi^\star_{ij}$ are $(n - 1)$-dimensional.
\par Let $C$ be a $\Pi_0$-conic. Then for $\pi_i \in C, \langle C \rangle \cap \tau_i = \bar \tau_i(C)$, or $\bar \tau_i$ if no confusion is possible, is the tangent space of $C$ at $\pi_i$. The projection of $C \setminus \lbrace \pi_0 \rbrace$ together with $\bar \tau_0 \cap \tau^\star_0$ is an $(n - 1)$-regulus $\mathcal{R}$. It is a family of maximal spaces of a Segre variety $\mathcal{S}_{1; n - 1}$. The space $\bar \tau^\star_0 = \bar \tau_0 \cap \tau^\star_0$ is the {\em space at infinity} of $\mathcal{R}$ (or $\mathcal{S}_{1; n - 1}$).
If $\pi^\star_i, \pi^\star_j \in \mathcal{R} \setminus \lbrace \bar \tau^\star_0 \rbrace, i \not = j$, then $\pi^\star_{ij} = \langle \pi^\star_i, \pi^\star_j \rangle \cap \tau^\star_0 = \bar \tau^\star_0$.
\par Let $\mathcal{R}_0, \mathcal{R}_1, \ldots$, respectively $\mathcal{S}^{(o)}_{1; n - 1}, \mathcal{S}^{(1)}_{1; n - 1}, \ldots$, be the
$(n - 1)$-reguli, respectively Segre varieties, defined by the pairs of elements of $\lbrace \pi^\star_1, \pi^\star_2, \ldots, \pi^\star_{q^{2n}} \rbrace$. Further, let $\xi_1, \xi_2, \dots$ be the $(q^n - 1)/(q - 1)$ $2n$-dimensional subspaces of $\Phi$ which contain $\tau^\star_0$.
\par Then $\xi_i \cap \pi^\star_k, k \in \lbrace 1, 2, \ldots, q^{2n} \rbrace$, is a point and $\xi_i \cap \mathcal{S}^{(j)}_{1; n - 1}$ consists of the space at infinity $\langle \mathcal{S}^{(j)}_{1; n - 1} \rangle \cap \tau^\star_0$ of $\mathcal{S}^{(j)}_{1; n - 1}$ and a maximal line of $\mathcal{S}^{(j)}_{1; n - 1}$ (distinct from the $\mathcal{R}_j$'s if $n = 2$), that is, a line intersecting the $q + 1$ elements of $\mathcal{R}_j$.
\par Let the affine space $\xi_i \setminus \tau^\star_0$ be denoted by $\AG^{(i)}(2n, q), i = 1, 2, \dots, (q^n - 1)/(q - 1)$. Let $\alpha_i: x \mapsto \pi^\star$, with $x$ any point of $\AG^{(i)}(2n, q)$ and $\pi^\star$ the element of $\lbrace \pi^\star_1, \pi^\star_2, \ldots, \pi^\star_{q^{2n}} \rbrace$ which contains $x$. Then on $\lbrace \pi^\star_1, \pi^\star_2, \ldots, \pi^\star_{q^{2n}} \rbrace$ a $2n$-dimensional affine space $\bar \AG^{(i)}(2n, q)$ is induced having as lines the $(n - 1)$-reguli $\mathcal{R}_k$ minus their space at infinity. For $q \not = 2$ we will prove that $\bar \AG^{(i)}(2n, q)$ is independent from $i$, which will appear not to be obvious for $q = 2$.
\par The $3n$-dimensional space $\langle \pi_0, \xi_i \rangle$, which contains $\tau_0$, will be denoted by $\bar \xi_i$, with $i = 1, 2, \ldots$. Each space $\pi_j \in O \setminus \lbrace \pi_0 \rbrace$, intersects $\bar \xi_i$ in exactly one point $z_j$. Let $Z_i$ be the set of all points $z_j$; then $\vert Z_i \vert = q^{2n}$.

\begin{lemma}
\label{lem6.15}
Any two points of $Z_i$ are contained in exactly one $\pi_0$-conic $\mathcal{C}$. Also, $\mathcal{C} \setminus \pi_0 \subset Z_i$ and $(\mathcal{C} \setminus \pi_0)^\delta$ is an affine line of $\AG^{(i)}(2n, q)$. In such a way all affine lines of $\AG^{(i)}(2n, q)$ are obtained.
\end{lemma}
{\em Proof} \quad
Let $x, y$ be two distinct points of $Z_i$. Let $\mathcal{C}$ be the unique $\pi_0$-conic which contains $x$ and intersects the space $\pi_j \in \O \setminus \lbrace \pi_0 \rbrace$ containing $y$. Let $\mathcal{C} \cap \tau_0 = \lbrace u \rbrace$. The tangent line of $\mathcal{C}$ at $u$ is contained in $\tau_0 \subset Z_i$, so $\langle \mathcal{C} \rangle$ is contained in $Z_i$, so $\mathcal{C} \cap \pi_j \in Z_i$, so $y \in \mathcal{C}$. Hence $x$ and $y$ are contained in a unique $\pi_0$-conic $\mathcal{C}$, and $\mathcal{C} \setminus \pi_0$ belongs to $Z_i$. Clearly $(\mathcal{C} \setminus \pi_0)^\delta$ is an affine line of $\AG^{(i)}(2n, q)$ and each affine line of $\AG^{(i)}(2n, q)$ is obtained in such a way. \eop \\

\par Consider distinct affine spaces $\AG^{(i)}(2n, q)$ and $\AG^{(j)}(2n, q)$. Let $\alpha_{ij}$ be the bijection from $\AG^{(i)}(2n, q)$ onto $\AG^{(j)}(2n, q)$ which maps the point $\xi_i \cap \pi^\star_k$ onto the point $\xi_j \cap \pi^\star_k$, with $k = 1, 2, \ldots, q^{2n}$.

\begin{lemma}
\label{lem6.16}
For $q \not = 2$ the bijection $\alpha_{ij}$ is a linear projectivity from $\AG^{(i)}(2n, q)$ onto $\AG^{(j)}(2n, q)$.
\end{lemma}
{\em Proof} \quad
If $L$ is a line of $\AG^{(i)}(2n, q)$, then $L$ and $L^{\alpha_{ij}}$ are maximal lines of one of the Segre varieties $\mathcal{S}_{1; n - 1}$. As $q \not = 2$ the mapping $\alpha_{ij}$ is a projectivity from $\AG^{(i)}(2n, q)$ onto $\AG^{(j)}(2n, q)$.
\par Let $\alpha_i: x \mapsto \pi^\star$, with $x \in \AG^{(i)}(2n, q)$ and $\pi^\star$ the element of $\lbrace \pi^\star_1, \pi^\star_2, \ldots, \pi^\star_{q^{2n}} \rbrace$ which contains $x$. Then on $\lbrace \pi^\star_1, \pi^\star_2, \ldots, \pi^\star_{q^{2n}} \rbrace$ a $2n$-dimensional space $\bar \AG^{(i)}(2n, q)$ is induced having as lines the $(n - 1)$-reguli $\mathcal{R}_k$ minus their space at infinity. From the preceding paragraph also follows that $\bar \AG^{(i)}(2n, q) = \bar \AG^{(j)}(2n, q)$. As $\alpha_i\alpha^{-1}_j= \alpha_{ij}$ it follows that $\alpha_{ij}$ preserves the cross-ratio of any four points on any affine line of $\AG^{(i)}(2n, q)$. Hence $\alpha_{ij}$ is linear. \eop \\

\begin{corollary}
\label{cor6.17}
$\bar \AG^{(i)}(2n, q) = \bar \AG^{(j)}(2n, q) = \bar \AG(2n, q)$ and all $(n - 1)$-reguli $\mathcal{R}_k$ which define a set of parallel lines of $\bar \AG(2n, q)$ have a common space at infinity (in any $\AG^{(l)}(2n, q)$ the maximal lines of the corresponding Segre varieties $\mathcal{S}^{(k)}_{1; n - 1}$ intersect $\tau^\star_0$ in a common point). Hence the corresponding $\Pi_0$-conics on $O$ have a common tangent space at $\pi_0$.
\end{corollary}

\par Let us consider a $\Pi_0$-set $(\Pi, \Gamma)$, with $\Pi = \lbrace \pi_0, \pi_1, \ldots, \pi_{q^n} \rbrace$ and $\Gamma = \lbrace \gamma_0, \gamma_1, \ldots, \gamma_{q^n} \rbrace$. The set of all points $\xi_i \cap \pi^\star_k$, with $k = 1, 2, \ldots, q^n$, will be denoted by $\xi^\star_i$.

\begin{lemma}
\label{lem6.18}
For $q \not = 2$ the set $\xi^\star_i$ is an $n$-dimensional affine subspace  $\AG^{(i)}(n, q)$ of $\AG^{(i)}(2n, q)$, with $i = 1, 2, \ldots, (q^n - 1)/(q - 1)$.
\end{lemma}
{\em Proof}
Let $z^\star_u \in \pi^\star_u, z^\star_v \in \pi^\star_v$, with $u, v \in \lbrace 1, 2, \ldots, q^n \rbrace$ and $u \not = v$. With the line $z^\star_uz^\star_v$ corresponds a $\pi_0$-conic $\mathcal{C}$, which is contained in $\bar \xi_i$ and which intersects the $q  + 1$ elements of the $\Pi_0$-conic determined by $\pi_0, \pi_u, \pi_v$; see \ref{lem6.10} and \ref{lem6.15}. Hence the affine line 
$z^\star_uz^\star_v \setminus \tau^\star_0$ belongs to $\xi^\star_i$. As $q \not = 2$, It follows that $\xi^\star_i$ is an affine $n$-dimensional subspace $\AG^{(i)}(n, q)$ of $\AG^{(i)}(2n, q)$, with $i = 1, 2, \ldots, (q^n - 1)/(q - 1)$. \eop \\

\begin{remark}
\label{rem6.19}
Consider again the $\Pi_0$-set $(\Pi, \Gamma)$, with $\Pi = \lbrace \pi_0, \pi_1, \ldots, \pi_{q^n} \rbrace$  and $\Gamma = \lbrace \gamma_0, \gamma_1, \ldots, \gamma_{q^n} \rbrace$. Assume that $\gamma_0 \cap \tau^\star_0 = \lbrace x^\star \rbrace$. Then the (projective completions of the ) lines of $\xi^\star_i$ through $x^\star$ correspond to the $\pi_0$-conics $\mathcal{C}$ in $\bar \xi_i  = \langle \pi_0, \xi_i \rangle$ intersecting $q$ elements of $\langle \pi_1, \pi_2, \ldots. \pi_{q^n} \rbrace$ and having their tangent line at $\mathcal{C} \cap \pi_0$ in the space $\gamma_0$.
\end{remark}

\begin{corollary}
\label{cor6.20} 
By \ref{lem6.18} the tangents of the $\pi_0$-conics in $\bar \xi_i$ of a $\Pi_0$-set all belong to a $(2n - 1)$-dimensional space $\beta_i$ containing $\pi_0$.
\end{corollary}

\subsection{MAIN THEOREMS}
\begin{theorem}
\label{thm6.21}
For $q \not = 2$ the pseudo-ovoid $O$ is good at $\pi_0$.
\end{theorem}
{\em Proof} \quad
We use the notations of Section 6.5. Let $q \not = 2$.
\par We introduced the projectivity $\alpha_{ij}$ from $\AG^{(i)}(2n, q) = \xi_i \setminus \tau^\star_0$ onto $\AG^{(j)}(2n, q) = \xi_j \setminus \tau^\star_0$. With each $(n - 1)$-regulus $\mathcal{R}_k$ (or Segre variety $\mathcal{S}^{(k)}_{1; n - 1}$) there corresponds a line $L$ of $\AG^{(i)}(2n, q)$ and a line $L^\prime$ of $\AG^{(j)}(2n, q)$; also $L^{\alpha_{ij}} = L^\prime$.
\par With a set of parallel lines of $\AG^{(i)}(2n, q)$ (and $\AG^{(j)}(2n, q)$) there corresponds a set of $(n - 1)$-reguli $\mathcal{R}_k$ with a common space at infinity in $\tau^\star_0$, and also a set of $\Pi_0$-conics with a common tangent space at $\pi_0$.
\par With any $\Pi_0$-set $(\Pi, \Gamma)$ there corresponds an affine $n$-dimensional subspace $\xi^\star_i = \AG^{(i)}(n, q)$ of $\AG^{(i)}(2n, q)$. Clearly $\AG^{(i)}(n, q)^{\alpha_{ij}} = \AG^{(j)}(n, q)$. Now let $\Pi = \lbrace \pi_0, \pi_1, \ldots, \pi_{q^n} \rbrace$, $\Gamma = \lbrace \gamma_0, \gamma_1, \ldots, \gamma_{q^n} \rbrace$, and $\gamma_0 \cap \tau^\star_0 = \lbrace x^\star_0 \rbrace$. The $\pi_0$-conics $\mathcal{C}$ in $\bar \xi^\star _i = \langle \pi_0, \xi_0 \rangle$ intersecting $q$ elements of $\Pi \setminus \lbrace \pi_0 \rbrace$ and having their tangent line at $\mathcal{C} \cap \pi_0$ in $\gamma_0$, correspond to the lines of $\AG^{(i)}(n, q)$ having $x^\star_0$ as point at infinity. The $\Pi_0$-conics $C$ containing the conics  $\mathcal{C}$ have a common $(2n - 1)$-space $\eta_i$ as tangent space at $\pi_0$. Let $\eta_i \cap \tau^\star_0 = \eta^\star_i$. These $\Pi_0$-conics are denoted by $C^i_1, C^i_2, \ldots, C^i_{q^{n - 1}}$.
\par Consider the $n$-dimensional space $\widetilde \gamma_0$ through $\pi_0$ in $\eta_i$ and consider the $\Pi_0$-set $(\widetilde \Pi, \widetilde \Gamma)$ defined by the pair $\lbrace \widetilde \gamma_0, \pi_1 \rbrace$, with $\pi_1 \in C^i_1$. Let $\widetilde \gamma_0 \cap \tau^\star_0 = \lbrace \widetilde x^\star_0 \rbrace$, with $\widetilde x^\star_0 \in \eta^\star_i$. As $\widetilde \gamma_0 \subset \eta_i$ and $\eta_i$ is the tangent space of $C^i_1$ at $\pi_0$, all elements of $C^i_1$ are in $\widetilde \Pi$. The $\pi_0$-conic $\mathcal{C}_1$ in $C^i_1$ with tangent line at $\mathcal{C}_1 \cap \pi_0$ in $\widetilde \gamma_0$, is projected onto a line containing $\widetilde x^\star_0$. This line belongs to some space $\xi_j$, and its affine part is contained in the $n$-dimensional affine space $\widetilde \xi^\star_j$ defined by $(\widetilde \Pi, \widetilde \Gamma)$. Notice that $\widetilde x^\star_0 = (x^\star_0)^{\alpha_{ij}}$ and that all lines of $\xi_j \setminus \tau^\star_0$ on $\widetilde x^\star_0$ are images for $\alpha_{ij}$ of all lines in $\xi_i \setminus \tau^\star_0$ on $x^\star_0$.
\par With $(\widetilde \Pi, \widetilde \Gamma)$ correspond affine $n$-dimensional spaces $\widetilde \xi^\star_i$ (in $\xi_i$) and $\widetilde \xi^\star_j$ (in $\xi_j$). The mapping $\alpha_{ji}$ sends all lines of $\widetilde \xi^\star_j$ on $\widetilde x^\star_0$ to all lines of $\widetilde \xi^\star_i$ on $x^\star_0$. The $(n - 1)$-reguli defined by these lines, that is $ \widetilde {\mathcal{R}}_1 = \mathcal{R}_1$ (the $(n - 1)$-regulus defined by $\mathcal{C}_1$ = the $(n - 1)$-regulus defined by the $\pi_0$-conic of $C_1$ in $\bar \xi_i$), $\widetilde {\mathcal{R}}_2, \ldots, \widetilde {\mathcal{R}}_{q^{n - 1}}$, have at infinity the common space $\eta^\star_i$.
\par By \ref{lem6.13} the space $\gamma_0$ defines a partition $\mathcal{P}$ of $O\setminus \lbrace \pi_0 \rbrace$ into sets of size $q^n$.
\par Let $(\Pi_k, \Gamma_k)$, with $k = 1, 2, \ldots, q^n$, be the $q^n$ $\Pi_0$-sets defined by $\mathcal{P}$ (one of these $\Pi_0$-sets is $(\Pi, \Gamma))$. For each $(\Pi_k, \Gamma_k)$ there are $q^{n - 1}$ $\pi_0$-conics $\mathcal{C}$ in $\bar \xi_i$ intersecting $q$ elements of $\Pi_k \setminus \lbrace \pi_0 \rbrace$ and having their tangent line at $\mathcal{C} \cap \pi_0$ in $\gamma_0$. The corresponding $\Pi_0$-conics have a common tangent space $\eta_i$ at $\pi_0$ (the corresponding lines of $\xi^\star_i$ have $x^\star_0$ as point at infinity). Hence each $n$-dimensional space $\gamma_0$ through $\pi_0$ determines for each $i$ a unique tangent space $\eta_i$. Any $n$-dimensional  space $\widetilde \gamma_0$ of $\eta_i$ through $\pi_0$ also determines the tangent space $\eta_i$.
\par Assume that $\eta_i$ and $\eta^\prime_i$ are two distinct such tangent spaces at $\pi_0$ and that $\eta_i \cap \eta^\prime_i$ contains an $n$-dimensional subspace $\zeta$ containing $\pi_0$. By the foregoing paragraphs $\zeta$ defines (for each $i$) a unique such tangent space, and so $\eta_i = \eta^\prime_i$. Hence for each $i$ all such tangent spaces corresponding to $\Pi_0$-sets have two by two just $\pi_0$ in common. 
\par So for each $i$ there arise $(q^n - 1)/(q - 1)$ $(2n - 1)$-dimensional subspaces $\eta_i$ containing $\pi_0$ having two by two just $\pi_0$ in common, and $(q^n - 1)/(q - 1)$ $(n - 1)$-dimensional subspaces $\eta^\star_i$ which form a partition of $\tau^\star_0$. The set of all $\eta_i$'s will be denoted by $\Sigma_i$, the partition containing the $\eta^\star_i$'s by $\Sigma^\star_i$.
\par Let $x^\star_0 \in \tau^\star_0$ and let $\eta^\star_i \in \Sigma^\star_i$ contain $x^\star_0$. Then $(x^\star_0)^{\alpha_{ij}} \in \eta^\star_i$. Let $\eta^\star_j \in \Sigma^\star_j$ contain $(x^\star_0)^{\alpha_{ij}}$. Then $(x^\star_0)^{\alpha_{ij}\alpha_{ji}} = x^\star_0 \in \eta^\star_j$. The spaces $\eta^\star_i$ and $\eta^\star_j$ are the spaces at infinity of a common $(n - 1)$-regulus. Hence $\eta^\star_i = \eta^\star_j$. So $\Sigma^\star_i = \Sigma^\star_j = \Sigma^\star$.
\par Consider again $(\Pi, \Gamma)$, with $\Pi = \lbrace \pi_0, \pi_1, \ldots, \pi_{q^n} \rbrace, \Gamma = \lbrace \gamma_0, \gamma_1, \ldots, \gamma_{q^n} \rbrace$, and let $(\Pi, \Gamma)$ be determined by $\lbrace \gamma_0, \pi_1 \rbrace$. The corresponding $\Pi_0$-conics, defined by the $\pi_0$-conics in $\bar \xi_i$, are $C^i_1, C^i_2, \ldots, C^i_{q^{n - 1}}$, with $C^i_1 = \lbrace \pi_0, \pi_1, \ldots, \pi_q \rbrace$. Let $C^\prime$ be the $\Pi_0$-conic in $\Pi$ containing $\pi_0, \pi_1, \pi_k$, with $k \in \lbrace q + 1, q + 2, \ldots, q^n \rbrace$. We have: $\gamma_0 \subset \tau_0, \gamma_0$ intersects $\langle \pi_1, \pi_k \rangle$ by \ref{lem6.10} , so $\gamma_0 \subset \langle \pi_0, \pi_1, \pi_k \rangle$, and hence the tangent space $\eta^{\prime \prime}$ of $C^\prime$ at $\pi_0$, that is $\langle \pi_0, \pi_1, \pi_k \rangle \cap \tau_0$, contains $\gamma_0$.
\par There is a $\pi_0$-conic $\mathcal{C}$ intersecting $\pi_0, \pi_1, \pi_k$ and having its tangent line at $\mathcal{C} \cap \pi_0$ in $\gamma_0$. Projecting $\mathcal{C}$ there arises a line $L$ of $\Phi$. Assume that $L$ is in $\xi_l$, so $\eta^{\prime \prime} \in \Sigma^\star_l = \Sigma^\star$. Hence $\eta^{\prime \prime}$ is tangent space at $\pi_0$ of each of the $\Pi_0$-conics defined by $\lbrace \gamma_0, \pi_1 \rbrace$. So $\eta^{\prime \prime}$ is the common tangent space at $\pi_0$ of the $\Pi_0$-conics defined by all triples $\pi_0, \pi_1, \pi_k$, and hence clearly of all $\Pi_0$-conics contained in $\Pi$. It follows that all $\Pi_0$-conics in $\Pi$ are contained in the $(3n - 1)$-dimensional space $\langle \eta^{\prime \prime}, \pi_0 \rangle$.
\par Each $\Pi_0$-conic $C^{\prime \prime}$ is contained in a $(\Pi, \Gamma)$-set. Hence the $(3n - 1)$-dimensional space $\langle C^{\prime \prime} \rangle$ contains $q^n + 1$ elements of $O$.
\par We conclude that $O$ is good at $\pi_0$. \eop \\

\begin{theorem}
\label{thm6.22}
Let $\mathcal{S} = T(n, 2n, q)$, with $q \not = 2$, be a translation generalized quadrangle with base point $(\infty)$, and having a regular line $L$ not incident wit $(\infty)$.
\begin{itemize}
\item[(i)] If $q$ is odd, then $\mathcal{S}$ is the point-line dual of the translation dual of a semifield flock translation generalized quadrangle.
\item[(ii)] If $q$ is even, then $\mathcal{S}$ is classical.
\end{itemize}
\end{theorem}
{\em Proof} \quad
Let $\mathcal{S} = T(O)$, with $O = \lbrace \pi_0, \pi_1, \ldots, \pi_{q{2n}} \rbrace$ and where the regular line $L$ is concurrent with $\pi_0$. By \ref {thm6.21} the pseudo-ovoid $O$ is good at $\pi_0$.
\par If $q$ is odd, then by \ref{cor4.3} the GQ $\mathcal{S}$ is the point-line dual of the translation dual of a semifield flock TGQ.
\par Let $q$ be even, let $T(O)$ be contained in $\PG(4n, q)$ and let $\langle O \rangle = \PG(4n - 1, q)$. Further, let $\langle \pi_0, \pi_1, \pi_2 \rangle = \PG(3n - 1, q)$ and let $\langle L, \PG(3n - 1, q) \rangle = \PG(3n, q)$. In $\PG(3n, q)$ a subquadrangle $\mathcal{S}^\prime = T(O^\prime)$ of order $q^n$ is induced, with $O^\prime$ the pseudo-oval consisting of the $q^n + 1$ elements of $O$ in $\PG(3n - 1, q)$.The line $L$ of $\mathcal{S}^\prime$ is regular, so by the action of the translation group of the TGQ $T(O^\prime)$ (with base point $(\infty)$) on $L$ we see that all lines of $\mathcal{S}^\prime$ concurrent with $\pi_0$ are regular. Now by 1.3.6(iv) of \cite {P&JAT:09}, all lines of $\mathcal{S}^\prime$ are regular. Hence, by the well-known theorem of Benson, see 5.2.1 of \cite {P&JAT:09}, the GQ $\mathcal{S}^\prime$ is isomorphic to the classical GQ arising from a nonsingular quadric of $\PG(4, q^n)$. Now by Brown and Lavrauw \cite {B&L:05} it follows that $\mathcal{S}$ is isomorphic to the classical GQ arising from a nonsingular elliptic quadric of $\PG(5, q^n)$; see also 7.8 of \cite {TTVM:06} \eop \\

\section{Applications}
\begin{theorem}
\label{thm7.1}
Let $\mathcal{S}$ be a GQ of order $(s, s^2), s \not = 1$, having two distinct translation points. Assume that for at least one of the corresponding TGQ's the kernel $\mathbb{K}$ has size at least 3.
\begin{itemize}
\item[(i)] If $s$ is odd, then $\mathcal{S}$ is the point-line dual of the translation dual of a semifield flock TGQ.
\item[(ii)] If $s$ is even, then $\mathcal{S}$ is classical.
\end{itemize}
\end{theorem}
{\em Proof} \quad
First assume that the two translation points $(\infty)$ and $(\infty)^\prime$ are collinear. Then all lines incident with $(\infty)^\prime$ are regular, so we can apply \ref{thm6.22} to the TGQ $\mathcal{S}^{(\infty)}$. 
\par Next let $(\infty)$ and $(\infty)^\prime$ be not collinear. There is a translation about $(\infty)$ which maps $(\infty)^\prime$ onto a point $(\infty)^{\prime \prime} \not = (\infty)^\prime$ collinear with $(\infty)^\prime$. So $\mathcal{S}$ has distinct collinear translation points, and the first part of the proof can be applied. \eop \\

\begin{remark}
\label{rem7.2}
Result \ref{thm7.1}, without the restriction $\vert \mathbb{K} \vert > 2$, was first proved by K. Thas \cite {KT:01}, relying on the classification of finite groups with a split $BN$-pair of rank 1 due independently to Hering, Kantor and Seitz \cite {HKS:72} and Shult \cite {S:72}.
\end{remark}

\begin{theorem}
\label{thm7.3} 
Let $\mathcal{S} = (P, B, \mathbf{I})$ be a GQ of order $(s, s^2)$, with $s \not = 1$, for which each point is a translation point. If no one of the TGQ's $\mathcal{S}^{(p)}$, with $p \in P$, has a kernel with 2 elements, then $\mathcal{S}$ is classical.
\end{theorem}
{\em Proof} \quad
It is clear that all lines of $\mathcal{S}$ are regular.
\par Let $\mathcal{S} = T(O)$. By \ref{thm6.21} the pseudo-ovoid $O$ is good at each of its elements. By 8.7.4(iii) of \cite {P&JAT:09}, the GQ $\mathcal{S}$ is isomorphic to a $T_3(\bar O)$ of Tits, with $\bar O$ some ovoid of $\PG(3, s)$. As all lines of $\mathcal{S}$ are regular, by 3.3.3(iii) of \cite {P&JAT:09}, the GQ is classical. \eop \\

\begin{remark}
\label{rem7.4}
So for $\vert \mathbb{K} \vert \not = 2$ this is the "missing part" in the monograph of Payne and Thas \cite {P&JAT:84,P&JAT:09}, so that now, for $\vert \mathbb{K} \vert \not = 2$, there is a pure geometrical proof of one of the main cases (the most complicated one) in the theorem of Fong and Seitz \cite {F&S:73,F&S:74} on $BN$-pairs of rank 2.
\end{remark} 

\begin{theorem}
\label{thm7.5}
Let $\mathcal{S} = (P, B, \mathbf{I})$ be a TGQ of order $(s, t)$, with $s \not = 1 \not = t$, and $s \not = t$, having $(\infty)$ as base point. If $\mathcal{S}$ has a regular line $L$ not incident with $(\infty)$, then $t = s^2$. Hence if the kernel has more than two elements, then \ref{thm6.22} can be applied.
\end{theorem}
{\em Proof} \quad
By 3.2 we know that for $s$ or $t$ even we necessarily have $t = s^2$. So from now on assume that $s$ and $t$ are odd. Let $L$ be concurrent with the line $\pi_0 \in O$ of $\mathcal{S} = T(O) = T(n, m, q)$. Similar to the case $t = s^2$ we prove that any $(3n - 1)$-dimensional space containing $\pi_0$ and two distinct elements $\pi_i, \pi_j \in O \setminus \lbrace \pi_0 \rbrace$, contains exactly $q^n + 1$ elements of $O$. Notice that $\langle \pi_0, \pi_i, \pi_j \rangle$ intersects no one of the remaining $q^m - q^n$ elements of $O$. So the spaces $\langle \pi_0, \pi_i, \pi_j \rangle$, with $\pi_j \in O \setminus \lbrace \pi_0, \pi_i \rbrace$ and $\pi_i$ fixed, determine a partition of $O \setminus \lbrace \pi_0, \pi_i \rbrace$ into sets of size $q^n - 1$. Hence $q^n - 1$ divides $q^m - 1$. It follows that $n$ divides $m$, and as $n \not = m$ we have $m \geq 2n$. As $m \leq 2n$ by Section 2.2, there follows that $m = 2n$, and so $t = s^2$. \eop \\

\section{the case $q = 2$}
All results in this paper also hold for $q = 2$ if \ref{lem6.16} can be proved for this case. In the proof of \ref{lem6.18} we also rely on $q \not = 2$, but for $q = 2$ this lemma can be proved relying on \ref{lem6.16}.

\end{document}